\documentclass[a4paper,12pt]{article}
\usepackage{amsfonts,amssymb,graphicx}

\newtheorem{lemma}{Lemma}

\newtheorem{remark}{Remark}

\begin{document}

\author{Fahd Jarad\\
{\small Department of Mathematics {\&} Computer Science,}\\
{\small \c{C}ankaya University, \"{O}gretmenler Cad. 14 06530, Balgat-Ankara, Turkey}\\
{\small e-mail: fahd@cankaya.edu.tr}\\
{\small and}\\
Octavian G. Mustafa\\
{\small Faculty of Mathematics, D.A.L., University of Craiova, Romania}\\
{\small e-mail: octaviangenghiz@yahoo.com}\\
{\small and}\\
Donal O'Regan\\
{\small Mathematics Department, National University of Ireland, Galway, Ireland}\\
{\small e-mail: donal.oregan@nuigalway.ie}\\
}

\title{Positive solutions of some elliptic differential equations with oscillating nonlinearity}
\date{}
\maketitle

\noindent\textbf{Abstract} We discuss the occurrence of positive solutions which decay to $0$ as $\vert x\vert\rightarrow+\infty$ to the differential equation $\Delta u+f(x,u)+g(\vert x\vert)x\cdot\nabla u=0$, $\vert x\vert>R>0$, $x\in\mathbb{R}^{n}$, where $n\geq 3$, $g$ is nonnegative valued and $f$ has alternating sign, by means of the comparison method. Our results complement several recent contributions from [M. Ehrnstr\"{o}m, O.G. Mustafa, On positive solutions of a class of nonlinear elliptic equations, Nonlinear Anal. TMA 67 (2007), 1147--1154]. 



\section{Introduction}

Consider the semilinear partial differential equation of second order
\begin{eqnarray}
\Delta u+f(x,u)+g(\vert x\vert)x\cdot\nabla u=0,\qquad x\in G_{R},\label{pde_1}
\end{eqnarray}
where $G_{R}=\{x\in\mathbb{R}^{n}:\vert x\vert>R\}$, $n\geq3$, $R>0$ and the functions $f,g$ verify the following smoothness assumptions: given $\alpha\in(0,1)$, $f\in C^{\alpha}(M\times J,\mathbb{R})$ for every compact set $M\subset G_{R}$ and every compact interval $J\subset\mathbb{R}$, and~$g\in C^{1}([R,+\infty),\mathbb{R})$.

In \cite{Constantin1996,Constantin1997}, A.~Constantin has established by means of comparison method \cite{GilbargTrudinger} --- see also \cite{NS} for a specifically designed approach --- that the equation (\ref{pde_1}) possesses a positive, decaying to $0$ as $\vert x\vert\rightarrow+\infty$, solution $u$ defined in $G_{A}$ for some $A\geq R$.

The main difficulty of the investigation from \cite{Constantin1996,Constantin1997} consists of the construction of the positive supersolution for the equation (\ref{pde_1}) since it requires solving an infinite-interval boundary value problem for a nonlinear ordinary differential equation.

It was also assumed there that
\begin{eqnarray}
0\leq f(x,u)\leq a(\vert x\vert)u,\qquad u\in[0,\varepsilon],\label{hyp_cu_semn}
\end{eqnarray}
for some $\varepsilon>0$, where $a:[R,+\infty)\rightarrow[0,+\infty)$ is continuous and such that
\begin{eqnarray}
\int_{R}^{+\infty}ra(r)dr<+\infty.\label{hyp_cu_semn_2}
\end{eqnarray}

The sign condition emerging from (\ref{hyp_cu_semn}), namely $uf(x,u)\geq0$, is essential for applying the maximum principle in the proofs from \cite{Constantin1996,Constantin1997}.

In \cite{E}, M. Ehrnstr\"{o}m was able to significantly improve the conclusions of \cite{Constantin1996,Constantin1997} in the case when $g(r)\geq0$ for all $r\geq R$ by noticing a special feature of the supersolutions to (\ref{pde_1}) which is the consequence of the particular form of the ordinary differential equation used for constructing the supersolutions:
\begin{eqnarray}
h^{\prime\prime}+p(s)\left(h^{\prime}-\frac{h}{s}\right)+q(s,h)=0,\qquad s\geq s_{0}>0.\label{ode_0}
\end{eqnarray}
A simplification of the proofs from \cite{Constantin1997,E} can be read in \cite{Mustafa2008}.

A further development of the techniques from \cite{Constantin1996,Constantin1997,E} has been done in \cite{EhrnstromMustafa}, where the comparison equation (\ref{ode_0}) was regarded as a {\it small perturbation} of the integrable ordinary differential equation
\begin{eqnarray*}
h^{\prime\prime}+p(s)\left(h^{\prime}-\frac{h}{s}\right)=0,\qquad s\geq s_{0}>0.
\end{eqnarray*}

One of the major shortcomings of \cite[Theorem 1]{EhrnstromMustafa} is that, unless $g$ verifies a technical hypothesis, the coefficient $a$ from (\ref{hyp_cu_semn}) will have to obey the condition
\begin{eqnarray*}
\int_{R}^{+\infty}r^{n-1}a(r)dr<+\infty,
\end{eqnarray*}
which is much more restrictive than (\ref{hyp_cu_semn_2}).

Our aim in this note is to give a positive (partial) answer to the following problem that, to the best of our knowledge, is open: assuming that $g$ is nonnegative valued everywhere, both functions $f,g$ are as smooth as necessary for the comparison method to work and the hypothesis (\ref{hyp_cu_semn}) is replaced by the condition
\begin{eqnarray*}
a_{1}(\vert x\vert)\leq f(x,u)\leq a_{2}(\vert x\vert),\qquad (x,u)\in{\cal D}\subseteq G_{R}\times[0,\varepsilon],
\end{eqnarray*}
where the continuous $a_{i}$'s have alternating sign and
\begin{eqnarray}
\int_{R}^{+\infty}r\vert a_{i}(r)\vert dr=+\infty,\qquad i\in\{1,2\},
\end{eqnarray}
with (if any) additional restrictions upon $f,g$, can one produce a positive solution $u$ of equation (\ref{pde_1}) such that $\lim\limits_{\vert x\vert\rightarrow+\infty}u(x)=0$ ?

The approach presented here relies on building a pair of subsolutions and supersolutions to the equation (\ref{pde_1}) with the help of some positive and bounded solutions of two differential equations of type (\ref{ode_0}). The proof of existence for a positive, vanishing at infinity, solution $u$ of equation (\ref{pde_1}) lying in between these subsolutions and supersolutions will follow then by a standard application of the classical comparison method.

This note consists of four sections. The second and third sections, of independent interest, deal with the problem of bounded solutions for the ordinary differential equation (\ref{ode_0}). The construction from the third section will be used for producing the needed family of functions $f,g$ in the last section.

\section{Positive solutions of certain equations (\ref{ode_0})}

Let $p:[s_{0},+\infty)\rightarrow[0,+\infty)$ and $q:[s_{0},+\infty)\rightarrow\mathbb{R}$ be two continuous functions such that $p$ is $L^{1}$ in $(s_{0},+\infty)$ and $q$ has alternating sign.

Suppose that the quantity
\begin{eqnarray}
z(s)=-\exp\left(\int_{s_0}^{s}p(\tau)d\tau\right)\int_{s_0}^{s}q(\tau)\exp\left(\int_{s_0}^{\tau}p(\xi)d\xi\right)d\tau\label{def_z}
\end{eqnarray}
is bounded in $[s_{0},+\infty)$. Then, the function $h:[s_{0},+\infty)\rightarrow\mathbb{R}$ with the formula
\begin{eqnarray}
h(s)=-s\int_{s}^{+\infty}\frac{z(\tau)}{\tau^2}d\tau\label{def_h}
\end{eqnarray}
will be a solution of the differential equation
\begin{eqnarray}
h^{\prime\prime}+p(s)\left(h^{\prime}-\frac{h}{s}\right)+\frac{q(s)}{s}=0,\qquad s\geq s_{0}.\label{ode_00}
\end{eqnarray}

Several restrictions will be imposed next on $p$, $q$ in order to obtain a positive solution $h$ in (\ref{def_h}).

\begin{lemma}\label{lem}
Set $\lambda=\int_{s_0}^{+\infty}p(\tau)d\tau<1$. Assume that there exist an increasing, unbounded from above sequence $(a_{m})_{m\geq2}$ of numbers from $[s_{0},+\infty)$ and another sequence $(\varepsilon_{m})_{m\geq1}$ of numbers from $(0,+\infty)$ such that
\begin{eqnarray}
q(s)>0,\thinspace s\in(a_{2m},a_{2m+1}),\qquad q(s)<0,\thinspace s\in(a_{2m+1},a_{2m+2}),\label{ode_hyp_0}
\end{eqnarray}
\begin{eqnarray}
\int_{a_{2m}}^{a_{2m+1}}q(s)ds\geq\left(1+3\int_{a_{2m}}^{+\infty}p(\tau)d\tau\right)\int_{a_{2m+1}}^{a_{2m+2}}\vert q(s)\vert ds\label{ode_hyp_1}
\end{eqnarray}
and
\begin{eqnarray}
\int_{a_{2m}}^{a_{2m+1}}q(s)ds\leq\varepsilon_{m}+\int_{a_{2m+1}}^{a_{2m+2}}\vert q(s)\vert ds\label{ode_hyp_2}
\end{eqnarray}
for all $m\geq 1$. 

Suppose also that $\sum\limits_{m=1}^{+\infty}\varepsilon_{m}=\varepsilon<+\infty$ and that
\begin{eqnarray}
\int_{a_{2m}}^{a_{2m+1}}q(s)ds\leq\delta<+\infty,\qquad m\geq 1.\label{ode_hyp_3}
\end{eqnarray}

Then, the function $z$ from (\ref{def_z}) is negative valued and bounded in the interval $[s_{0},+\infty)$ and, consequently, $h$ is a positive valued and bounded solution of equation (\ref{ode_00}) with the quantity $\frac{h(s)}{s}$ decreasing to $0$ as $s\rightarrow+\infty$.
\end{lemma}

{\bf Proof.} We start by noticing that $\exp x\leq1+3x$  for all $x\in[0,1]$. Since $\lambda\in[0,1)$, we have
\begin{eqnarray*}
\exp\left(\int_{a_{2m}}^{a_{2m+2}}p(\tau)d\tau\right)\leq1+3\int_{a_{2m}}^{a_{2m+2}}p(\tau)d\tau,\quad m\geq1.
\end{eqnarray*}

Taking $a_{2}=s_{0}$ for simplicity, we have the estimate
\begin{eqnarray*}
&&\int_{a_{2m}}^{a_{2m+2}}q(s)\exp\left(\int_{s_0}^{s}p(\tau)d\tau\right)ds\\
&&=\left(\int_{a_{2m}}^{a_{2m+1}}+\int_{a_{2m+1}}^{a_{2m+2}}\right)q(s)\exp\left(\int_{s_0}^{s}p(\tau)d\tau\right)ds\\
&&\geq\int_{a_{2m}}^{a_{2m+1}}q(s)ds\cdot\exp\left(\int_{s_0}^{a_{2m}}p(\tau)d\tau\right)\\
&&-\int_{a_{2m+1}}^{a_{2m+2}}\vert q(s)\vert ds\cdot\exp\left(\int_{s_0}^{a_{2m+2}}p(\tau)d\tau\right)\\
&&\geq\exp\left(\int_{s_0}^{a_{2m}}p(\tau)d\tau\right)\left[\int_{a_{2m}}^{a_{2m+1}}q(s)ds\right.\\
&&\left.-\left(1+3\int_{a_{2m}}^{a_{2m+2}}p(\tau)d\tau\right)\int_{a_{2m+1}}^{a_{2m+2}}\vert q(s)\vert ds\right]\\
&&\geq0,\qquad m\geq1,
\end{eqnarray*}
by taking into account the inequality (\ref{ode_hyp_1}).

Thus we have
\begin{eqnarray*}
\int_{s_0}^{a_{2m}}q(s)\exp\left(\int_{s_0}^{s}p(\tau)d\tau\right)ds&=&\sum\limits_{k=1}^{m-1}\int_{a_{2k}}^{a_{2k+2}}q(s)\exp\left(\int_{s_0}^{s}p(\tau)d\tau\right)ds\\
&\geq&0,\qquad m\geq2.
\end{eqnarray*}

As the local maxima of $z(s)$ are attained when $s=a_{2m}$, the preceding computations establish that $z$ is negative valued in $[s_{0},+\infty)$.

Further, we have the estimate
\begin{eqnarray*}
&&\int_{a_{2m}}^{a_{2m+2}}q(s)\exp\left(\int_{s_0}^{s}p(\tau)d\tau\right)ds\\
&&=\left(\int_{a_{2m}}^{a_{2m+1}}+\int_{a_{2m+1}}^{a_{2m+2}}\right)q(s)\exp\left(\int_{s_0}^{s}p(\tau)d\tau\right)ds\\
&&\leq\exp\left(\int_{s_0}^{a_{2m+1}}p(\tau)d\tau\right)\left(\int_{a_{2m}}^{a_{2m+1}}q(s)ds-\int_{a_{2m+1}}^{a_{2m+2}}\vert q(s)\vert ds\right)\\
&&\leq(\exp\lambda)\varepsilon_{m},
\end{eqnarray*}
according to (\ref{ode_hyp_2}), and respectively (recall (\ref{ode_hyp_3}))
\begin{eqnarray*}
\vert z(s)\vert&\leq&\sum\limits_{k=1}^{m-1}\int_{a_{2k}}^{a_{2k+2}}q(\tau)\exp\left(\int_{s_0}^{\tau}p(\xi)d\xi\right)d\tau\\
&+&\int_{a_{2m}}^{a_{2m+1}}q(\tau)\exp\left(\int_{s_0}^{\tau}p(\xi)d\xi\right)d\tau\\
&<&(\varepsilon+\delta)\exp\lambda,\qquad s\in[a_{2m},a_{2m+2}].
\end{eqnarray*}

The proof is complete. $\square$

We emphasize two particular cases of Lemma \ref{lem}. 

The first one is when the continuous function $q$ is $L^{1}$ in $(s_{0},+\infty)$ and satisfies the conditions (\ref{ode_hyp_0}), (\ref{ode_hyp_1}). Then, we can take
\begin{eqnarray*}
\varepsilon_{m}=\int_{a_{2m}}^{a_{2m+1}} q(s) ds,\quad m\geq1.
\end{eqnarray*}
Here,
\begin{eqnarray*}
\varepsilon=\sum\limits_{m=1}^{+\infty}\varepsilon_{m}\leq\sum\limits_{m=1}^{+\infty}\int_{a_{2m}}^{a_{2m+2}} \vert q(s)\vert ds=\Vert q\Vert_{L^{1}}
\end{eqnarray*}
and we may use $\delta=\varepsilon$.

In the second case, $q$ satisfies the conditions (\ref{ode_hyp_0}), (\ref{ode_hyp_1}), (\ref{ode_hyp_2}) together with
\begin{eqnarray}
0<q_{-}\leq\int_{a_{2m+1}}^{a_{2m+2}}\vert q(s)\vert ds<q_{+}<+\infty,\qquad m\geq1.\label{hyp_q_pde_0}
\end{eqnarray}
Here, $\int_{s_0}^{+\infty}\vert q(s)\vert ds\geq\sum\limits_{m=1}^{+\infty} q_{-}=+\infty$ and we may use $\delta=\varepsilon+q_{+}$. 

\begin{remark}
We deduce from (\ref{ode_hyp_1}), (\ref{ode_hyp_2}) that
\begin{eqnarray*}
\int_{a_{2m}}^{+\infty}p(\tau)d\tau\cdot\int_{a_{2m+1}}^{a_{2m+2}}\vert q(s)\vert ds\leq\varepsilon_{m},\qquad m\geq1.
\end{eqnarray*}

This means that, when (\ref{hyp_q_pde_0}) holds, an additional hypothesis must be introduced for the function $p$, namely
\begin{eqnarray*}
\sum\limits_{m=1}^{+\infty}\int_{a_{2m}}^{+\infty}p(\tau)d\tau\leq\frac{1}{q_{-}}\sum\limits_{m=1}^{+\infty}\varepsilon_{m}=\frac{\varepsilon}{q_{-}}<+\infty.
\end{eqnarray*}

To give a computational particular case of this situation, assume that
\begin{eqnarray}
\inf\limits_{m\geq1}(a_{2m+2}-a_{2m})=A>0\label{hyp_q_pde_00}
\end{eqnarray}
and
\begin{eqnarray}
\int_{s_0}^{+\infty}(s-s_{0})p(s)ds=\int_{s_0}^{+\infty}\int_{s}^{+\infty}p(\tau)d\tau ds<+\infty.\label{hyp_q_pde_1}
\end{eqnarray}
Then, we have the estimate
\begin{eqnarray*}
\int_{s_0}^{+\infty}(s-s_{0})p(s)ds&=&\sum\limits_{m=2}^{+\infty}\int_{a_{2m-2}}^{a_{2m}}\int_{s}^{+\infty}p(\tau)d\tau ds\\
&\geq&\sum\limits_{m=2}^{+\infty}\int_{a_{2m-2}}^{a_{2m}}\int_{a_{2m}}^{+\infty}p(\tau)d\tau ds\\
&&\geq A\sum\limits_{m=2}^{+\infty}\int_{a_{2m}}^{+\infty}p(\tau)d\tau ds.
\end{eqnarray*}
\end{remark}

\section{Example}

Set $a_{m}=m\pi$, where $m\geq1$. Thus, recalling (\ref{hyp_q_pde_00}), we have $A=2\pi>0$.

Assume also that $\lambda=\int_{s_0}^{+\infty}p(\tau)d\tau<1$ and (\ref{hyp_q_pde_1}) holds.

Set
\begin{eqnarray}
\left\{
\begin{array}{ll}
0<q_{-}<q_{+}<+\infty,\\
6\leq\gamma<\sigma<+\infty,\\
0\leq\eta<\theta<+\infty
\end{array}
\right.
\label{ex_const}
\end{eqnarray}
and introduce the sequences $(c_{m})_{m\geq1}$, $(d_{m})_{m\geq1}$ via the restrictions
\begin{eqnarray*}
q_{-}\leq\frac{\pi}{2}d_{m}\leq q_{+}
\end{eqnarray*}
and
\begin{eqnarray*}
&&d_{m}+\gamma\cdot\frac{q_{+}}{\pi}\int_{a_{2m}}^{+\infty}p(\tau)d\tau+\eta\cdot\frac{2^{1-m}}{\pi}\leq c_{m}\\
&&\leq d_{m}+\sigma\cdot\frac{q_{+}}{\pi}\int_{a_{2m}}^{+\infty}p(\tau)d\tau+\theta\cdot\frac{2^{1-m}}{\pi}
\end{eqnarray*}
for all $m\geq1$.

Introduce the function $q:[s_{0},+\infty)\rightarrow\mathbb{R}$ with the formula
\begin{eqnarray}
q(s)=
\left\{
\begin{array}{ll}
c_{m}\sin^{2}s,\thinspace s\in[a_{2m},a_{2m+1}],\\
-d_{m}\sin^{2}s,\thinspace s\in[a_{2m+1},a_{2m+2}],
\end{array}
\right.\label{ex_form_q}
\end{eqnarray}
and notice that it is continuously differentiable and bounded in $[s_{0},+\infty)$ from the estimate
\begin{eqnarray}
\vert q(s)\vert&\leq&\max\{c_{m},d_{m}\}=c_{m}\nonumber\\
&\leq&d_{m}+\sigma\cdot\frac{q_{+}}{\pi}\lambda+\theta\cdot\frac{1}{\pi}\nonumber\\
&\leq&\frac{1}{\pi}[(2+\sigma\lambda)q_{+}+\theta]<+\infty.\label{ex_marg_q}
\end{eqnarray}

We have also
\begin{eqnarray*}
\int_{a_{2m}}^{a_{2m+1}}q(s)ds=\frac{\pi}{2}c_{m},\quad\int_{a_{2m+1}}^{a_{2m+2}}\vert q(s)\vert ds=\frac{\pi}{2}d_{m}.
\end{eqnarray*}

Now,
\begin{eqnarray*}
&&\int_{a_{2m}}^{a_{2m+1}}q(s)ds-\int_{a_{2m+1}}^{a_{2m+2}}\vert q(s)\vert ds=\frac{\pi}{2}(c_{m}-d_{m})\\
&&\geq\frac{\gamma}{2}\cdot q_{+}\int_{a_{2m}}^{+\infty}p(\tau)d\tau+\eta\cdot 2^{-m}\geq\frac{\gamma}{2}\cdot q_{+}\int_{a_{2m}}^{+\infty}p(\tau)d\tau\\
&&\geq3\int_{a_{2m}}^{+\infty}p(\tau)d\tau\cdot\int_{a_{2m+1}}^{a_{2m+2}}\vert q(s)\vert ds
\end{eqnarray*}
which means that (\ref{ode_hyp_1}) holds.

Further,
\begin{eqnarray*}
&&\int_{a_{2m}}^{a_{2m+1}}q(s)ds-\int_{a_{2m+1}}^{a_{2m+2}}\vert q(s)\vert ds\\
&&\leq\frac{\sigma}{2}\cdot q_{+}\int_{a_{2m}}^{+\infty}p(\tau)d\tau+\theta\cdot 2^{-m}=\varepsilon_{m}
\end{eqnarray*}
which leads to (\ref{ode_hyp_2}).

In conclusion, the function $q$ from (\ref{ex_form_q}) fulfills all the requirements of Lemma \ref{lem} and generates a positive and bounded solution $h$ for the equation (\ref{ode_00}) such that
\begin{eqnarray*}
s\frac{d}{ds}\left(\frac{h(s)}{s}\right)=h^{\prime}(s)-\frac{h(s)}{s}=\frac{z(s)}{s}<0,\qquad s\geq s_{0}.
\end{eqnarray*}
For the significance of this estimate, see \cite{E,EhrnstromMustafa}.

Two important features of the example follow from the next computations, namely
\begin{eqnarray*}
&&\int_{s_0}^{+\infty}\frac{\vert q(s)\vert}{s}ds\geq\sum\limits_{m=1}^{+\infty}\int_{a_{2m+1}}^{a_{2m+2}}\frac{\vert q(s)\vert}{s}ds\\
&&\geq\sum\limits_{m=1}^{+\infty}\int_{a_{2m+1}+\frac{\pi}{4}}^{a_{2m+2}-\frac{\pi}{4}}d_{m}\cdot\frac{\sin^{2}s}{s}ds\geq\sum\limits_{m=1}^{+\infty}\int_{a_{2m+1}+\frac{\pi}{4}}^{a_{2m+2}-\frac{\pi}{4}}d_{m}\cdot\frac{1/2}{s}ds\\
&&\geq\sum\limits_{m=1}^{+\infty}d_{m}\cdot\frac{1/2}{(2m+2)\pi-\frac{\pi}{4}}\cdot\frac{\pi}{2}\\
&&\geq\sum\limits_{m=1}^{+\infty}d_{m}\cdot\frac{1}{8(m+1)}\geq\frac{2}{\pi}q_{-}\cdot(+\infty)=+\infty
\end{eqnarray*}
and, given $\varsigma>0$ and recalling (\ref{ex_marg_q}),
\begin{eqnarray*}
\int_{s_0}^{+\infty}\frac{\vert q(s)\vert}{s^{1+\varsigma}}ds\leq\int_{s_0}^{+\infty}\frac{\Vert q\Vert_{\infty}}{s^{1+\varsigma}}ds<+\infty.
\end{eqnarray*}

To deal with the problem stated in the introduction, we introduce now a pair $(q_1,q_2)$ of functions verifying (\ref{ex_const}), (\ref{ex_form_q}) such that
\begin{eqnarray}
q_{1}(s)\leq q_{2}(s),\qquad s\geq s_{0}.\label{ex_ineq_per_pde}
\end{eqnarray}

We shall use the upper index $i$ when referring to the constants from (\ref{ex_const}) that characterize the function $q_i$, where $i\in\{1,2\}$.

Set
\begin{eqnarray*}
q_{-}^{i}=q_{-},\quad q_{+}^{i}=q_{+},
\end{eqnarray*}
and
\begin{eqnarray*}
\sigma^{1}<\gamma^{2},\quad\theta^{1}<\eta^{2}.
\end{eqnarray*}

Set also $\alpha\in(0,\gamma^{2}-\sigma^{1})$ and $\beta\in(0,\eta^{2}-\theta^{1})$ small enough to have
\begin{eqnarray*}
q_{-}+\frac{\alpha}{2}\cdot q_{+}\lambda+\frac{\beta}{2}<q_{+}
\end{eqnarray*}
and notice that this restriction implies
\begin{eqnarray*}
\frac{2}{\pi}q_{-}+\alpha\cdot\frac{q_{+}}{\pi}\int_{a_{2m}}^{+\infty}p(\tau)d\tau+\beta\cdot\frac{2^{1-m}}{\pi}<\frac{2}{\pi}q_{+},\quad m\geq1.
\end{eqnarray*}

Introduce the sequences $(d_{m}^{i})_{m\geq1}$ via the restrictions
\begin{eqnarray*}
&&\frac{2}{\pi}q_{-}\leq d_{m}^{1}-\alpha\cdot\frac{q_{+}}{\pi}\int_{a_{2m}}^{+\infty}p(\tau)d\tau-\beta\cdot\frac{2^{1-m}}{\pi}\\
&&\leq d_{m}^{2}\leq d_{m}^{1}\leq\frac{2}{\pi}q_{+}.
\end{eqnarray*}

The sequences $(c_{m}^{i})_{m\geq1}$ are given by the inequalities
\begin{eqnarray*}
&&d_{m}^{1}+\gamma^{1}\cdot\frac{q_{+}}{\pi}\int_{a_{2m}}^{+\infty}p(\tau)d\tau+\eta^{1}\cdot\frac{2^{1-m}}{\pi}\leq c_{m}^{1}\\
&&\leq d_{m}^{1}+\sigma^{1}\cdot\frac{q_{+}}{\pi}\int_{a_{2m}}^{+\infty}p(\tau)d\tau+\theta^{1}\cdot\frac{2^{1-m}}{\pi}\\
&&\leq d_{m}^{2}+\alpha\cdot\frac{q_{+}}{\pi}\int_{a_{2m}}^{+\infty}p(\tau)d\tau+\beta\cdot\frac{2^{1-m}}{\pi}\\
&&+\sigma^{1}\cdot\frac{q_{+}}{\pi}\int_{a_{2m}}^{+\infty}p(\tau)d\tau+\theta^{1}\cdot\frac{2^{1-m}}{\pi}\\
&&\leq d_{m}^{2}+\gamma^{2}\cdot\frac{q_{+}}{\pi}\int_{a_{2m}}^{+\infty}p(\tau)d\tau+\eta^{2}\cdot\frac{2^{1-m}}{\pi}\leq c_{m}^{2}\\
&&\leq d_{m}^{2}+\sigma^{2}\cdot\frac{q_{+}}{\pi}\int_{a_{2m}}^{+\infty}p(\tau)d\tau+\theta^{2}\cdot\frac{2^{1-m}}{\pi}
\end{eqnarray*}
for all $m\geq1$. 

Notice the inequalities
\begin{eqnarray*}
c_{m}^{2}\geq c_{m}^{1}\geq d_{m}^{1}\geq d_{m}^{2}>0,\qquad m\geq1,
\end{eqnarray*}
which help establishing (\ref{ex_ineq_per_pde}).

\section{Positive solution of the equation (\ref{pde_1})}

Given $n\geq3$ and $R>0$, introduce the function $\beta:[s_{0},+\infty)\rightarrow[R,+\infty)$ with the formula
\begin{eqnarray*}
\beta(s)=\left(\frac{s}{n-2}\right)^{\frac{1}{n-2}},\qquad s_{0}>(n-2)R^{n-2}.
\end{eqnarray*}

Consider also the smooth functions $v,h$ connected by
\begin{eqnarray*}
v(x)=\frac{h(s)}{s},\qquad x\in\mathbb{R}^{n},\thinspace s\geq s_{0},\thinspace\vert x\vert=\beta(s).
\end{eqnarray*}
Obviously, $v$ is radially symmetric.

It can be established readily that
\begin{eqnarray}
&&\Delta v+f(x,v)+g(\vert x\vert)x\cdot\nabla v\label{transfo}\\
&&=\frac{n-2}{\beta(s)\beta^{\prime}(s)}\left[h^{\prime\prime}(s)+\beta(s)\beta^{\prime}(s)g(\beta(s))\left(h^{\prime}(s)-\frac{h(s)}{s}\right)\right.\nonumber\\
&&\left.+\frac{1}{n-2}\beta(s)\beta^{\prime}(s)f(x,v)\right],\nonumber
\end{eqnarray}
see also \cite[p. 1150]{EhrnstromMustafa}.

Assume now that the function $f$ verifies the following hypothesis
\begin{eqnarray}
a_{1}(\vert x\vert)\leq f(x,u)\leq a_{2}(\vert x\vert),\qquad (x,u)\in{\cal D},\label{ribbon}
\end{eqnarray}
where the continuous functions $a_{i}:[R,+\infty)\rightarrow\mathbb{R}$ have alternating sign and
\begin{eqnarray*}
{\cal D}=\left\{(x,u):x\in\mathbb{R}^{n},\thinspace\vert x\vert>R,\thinspace\frac{h_{1}(s)}{s}\leq u\leq\frac{h_{2}(s)}{s}\mbox{ for }\beta(s)=\vert x\vert\right\}
\end{eqnarray*}
and
\begin{eqnarray*}
\left\{
\begin{array}{ll}
h_{i}(s)=-s\int_{s}^{+\infty}\frac{z_{i}(\tau)}{\tau^{2}}d\tau,\\
\\ z_{i}(s)=-\exp\left(-\int_{s_0}^{s}p(\tau)d\tau\right)\int_{s_0}^{s}q_{i}(\tau)\exp\left(\int_{s_0}^{\tau}p(\xi)d\xi\right)d\tau,
\end{array}
\right.
\end{eqnarray*}
and
\begin{eqnarray*}
p(s)=\beta(s)\beta^{\prime}(s)g(\beta(s)),\qquad q_{i}(s)=\frac{s}{n-2}\beta(s)\beta^{\prime}(s)a_{i}(\beta(s))
\end{eqnarray*}
for $s\geq s_{0}$ and $i\in\{1,2\}$. Here, the functions $p$, $q_i$ are supposed to verify the requirements of Lemma \ref{lem}.

An example of such a pair $(q_1,q_2)$ has been given in the preceding section.

Notice that the functions $h_{i}$ are positive valued and bounded solutions of the linear differential equations
\begin{eqnarray*}
h^{\prime\prime}+p(s)\left(h^{\prime}-\frac{h}{s}\right)+\frac{q_{i}(s)}{s}=0,\qquad s\geq s_{0},
\end{eqnarray*}
which belong to the class (\ref{ode_0}).

Given the functions $v_i$ with $v_{i}(x)=\frac{h_{i}(s)}{s}$ for $\beta(s)=\vert x\vert>R$, we remark that via (\ref{transfo})
\begin{eqnarray*}
&&\Delta v_{1}+f(x,v_{1})+g(\vert x\vert)x\cdot\nabla v_{1}\\
&&\geq\frac{n-2}{\beta(s)\beta^{\prime}(s)}\left[h_{1}^{\prime\prime}(s)+\beta(s)\beta^{\prime}(s)g(\beta(s))\left(h_{1}^{\prime}(s)-\frac{h_{1}(s)}{s}\right)\right.\\
&&\left.+\frac{1}{n-2}\beta(s)\beta^{\prime}(s)a_{1}(\beta(s))\right]\\
&&=0
\end{eqnarray*}
and respectively that
\begin{eqnarray*}
&&\Delta v_{2}+f(x,v_{2})+g(\vert x\vert)x\cdot\nabla v_{2}\\
&&\leq\frac{n-2}{\beta(s)\beta^{\prime}(s)}\left[h_{2}^{\prime\prime}(s)+\beta(s)\beta^{\prime}(s)g(\beta(s))\left(h_{2}^{\prime}(s)-\frac{h_{2}(s)}{s}\right)\right.\\
&&\left.+\frac{1}{n-2}\beta(s)\beta^{\prime}(s)a_{2}(\beta(s))\right]\\
&&=0.
\end{eqnarray*}
Thus, $v_{1}$ is a subsolution and $v_{2}$ is a supersolution of the equation (\ref{pde_1}).

We have $v_{1}(x)\leq v_{2}(x)$ for $x\in G_{R}$. According to the classical comparison method \cite{GilbargTrudinger}, the equation (\ref{pde_1}) has a solution, non necessarily of radial symmetry, such that
\begin{eqnarray*}
v_{1}(x)\leq u(x)\leq v_{2}(x),\qquad x\in G_{R}.
\end{eqnarray*}
The functions $v_i$ are positive valued, so we conclude that this solution $u$ is also positive valued and behaves like
\begin{eqnarray*}
u(x)=O\left(s^{-1}\right)=O\left(\vert x\vert^{2-n}\right)
\end{eqnarray*}
when $\vert x\vert\rightarrow+\infty$.

In the end, let us recall the function $q$ from (\ref{ex_form_q}). The functions $a_i$ from the double inequality (\ref{ribbon}) are modeled  by
\begin{eqnarray*}
a(\vert x\vert)=a(\beta(s))=\frac{q(s)}{\frac{s}{n-2}\beta(s)\beta^{\prime}(s)}=\frac{q(s)}{[\beta(s)]^{n-1}\beta^{\prime}(s)}.
\end{eqnarray*}

The restrictions regarding $p,q$ from the example yield
\begin{eqnarray*}
\int_{s_0}^{+\infty}sp(s)ds<+\infty,\qquad\int_{s_0}^{+\infty}\frac{\vert q(s)\vert}{s} ds=+\infty
\end{eqnarray*}
and also
\begin{eqnarray*}
\int_{s_{0}}^{+\infty}\frac{\vert q(s)\vert}{s^{1+\varsigma}}ds<+\infty\qquad\mbox{when }\varsigma>0.
\end{eqnarray*}

Translated for $a_{i},g$, they read as
\begin{eqnarray*}
\int_{R}^{+\infty}r^{n-1}g(r)dr<+\infty
\end{eqnarray*}
and
\begin{eqnarray*}
\int_{R}^{+\infty}r\vert a_{i}(r)\vert dr=+\infty,\quad\int_{R}^{+\infty}r^{1-\varsigma(n-2)}\vert a_{i}(r)\vert dr<+\infty.
\end{eqnarray*}

As a result, the methods from \cite{Constantin1997}--\cite{EhrnstromMustafa} or \cite{Mustafa2008} are not applicable here.

\section{Acknowledgements.} This note was completed during the visit of OGM to The Erwin Schr\"{o}dinger International Institute for Mathematical Physics in October 2009. OGM is deeply indebted to Prof. Adrian Constantin from Vienna University and to the ESI staff for their support.

\end{document}